\documentclass[12pt]{amsart}
\usepackage{amscd,amssymb}
\usepackage[arrow,matrix]{xy}
\usepackage[plainpages,backref,urlcolor=blue]{hyperref}

\theoremstyle{plain}
\newtheorem{thm}[subsection]{Theorem}
\newtheorem{lem}[subsection]{Lemma}
\newtheorem{prop}[subsection]{Proposition}
\newtheorem{cor}[subsection]{Corollary}

\theoremstyle{definition}
\newtheorem{rk}[subsection]{Remark}
\newtheorem{definition}[subsection]{Definition}
\newtheorem{ex}[subsection]{Example}
\newtheorem{question}[subsection]{Question}

\numberwithin{equation}{section}
\newcommand{\OO}{{\mathcal O}}

\newcommand{\I}{{\mathcal I}}

\newcommand{\C}{\mathbb{C}}

\newcommand{\PP}{\mathbb{P}}
\newcommand{\AAA}{\mathbb{A}}

\DeclareMathOperator{\im}{im}

\DeclareMathOperator{\defect}{def}

\begin{document}
%\date{June 4, 2009}

\title [Syzygies of Jacobian ideals and weighted homogeneous singularities]
{Syzygies of Jacobian ideals and weighted homogeneous singularities }

\author[Alexandru Dimca]{Alexandru Dimca$^1$}
\address{Univ. Nice Sophia Antipolis, CNRS,  LJAD, UMR 7351, 06100 Nice, France. }
\email{dimca@unice.fr}

\author[Gabriel Sticlaru]{Gabriel Sticlaru}
\address{Faculty of Mathematics and Informatics,
Ovidius University,
Bd. Mamaia 124, 900527 Constanta,
Romania}
\email{gabrielsticlaru@yahoo.com }
\thanks{$^1$ Supported by Institut Universitaire de France}



\begin{abstract} Let $V$ be a projective hypersurface having only isolated  singularities. We show that these singularities are weighted homogeneous if and only if the Koszul syzygies among the partial derivatives of an equation for $V$ are exactly the syzygies with a generic first component vanishing on the singular locus subscheme of $V$. This yields in particular a positive answer in this setting to a question raised by Morihiko Saito and the first author.
Finally we explain how our result can be used to improve the listing of Jacobian syzygies of a given degree by a computer algebra system such as Singular, CoCoA or Macaulay2.

\end{abstract}

\maketitle

%\tableofcontents

\section{Introduction and statement of results} \label{sec:intro}

Let $X$ be the complex projective space $\PP^n$ and consider the associated graded polynomial algebra $S= \oplus _kS_k $ with $S_k=H^0(X, \OO_X(k))$. For a nonzero section $f \in S_N$ with $N>1$, we consider the hypersurface $V=V(f)$ in $X$ given by the zero locus of $f$and let $Y$ denote the singular locus of $V$, endowed with its natural scheme structure, see \cite{DBull}.
We assume in this paper that $V$ has isolated singularities.

Let $\I_Y \subset \OO_X$ be the ideal sheaf defining this subscheme $Y \subset X$ and consider the graded ideal $I= \oplus_k I_k$ in $S$ with $I_k=H^0(X, \I_Y(k))$. Let $Z=Spec(S)$ be the corresponding affine space $\AAA^{n+1}$ and denote by $\Omega^j=H^0(Z, \Omega_Z^k)$ the $S$-module of global, regular $k$-forms on $Z$. Using a linear coordinate system $x=(x_0,...,x_n)$ on $X$, one sees that there is a natural isomorphism of $S$-modules
\begin{equation} 
\label{eq1}
\Omega^j=S^{n+1 \choose j}
\end{equation}
which is used to put a grading (independent of the choice of $x$) on the modules $\Omega^j$, i.e.  a differential form 
\begin{equation} 
\label{eq2}
\omega=\sum _K \omega_K(x) dx_K
\end{equation}
 is homogeneous of degree $m$ if all the coefficients $\omega_K(x)$ are in $ S_m$ for all multi-indices $K=(k_1<...<k_j)$ and $dx_K=dx_{k_1} \wedge ... \wedge dx_{k_j} $.

Since $f$ can be thought of as a homogeneous polynomial of degree $N$ on $Z$, it follows that there is a well defined differential $1$-form $df \in \Omega^1$. Using this we define two graded $S$-submodules in $\Omega^n$, namely
\begin{equation} 
\label{eqAR}
AR(f)=\ker \{df \wedge : \Omega^n \to \Omega^{n+1}\}
\end{equation}
and
\begin{equation} 
\label{eqKR}
KR(f)=\im \{df \wedge : \Omega^{n-1}\to   \Omega^n\}.
\end{equation}
If one computes in a coordinate system $x$, then $AR(f)$ is the module of {\it all relations} of the type 
\begin{equation} 
\label{eqR}
R: a_0f_{x_0}+...a_nf_{x_n}=0,
\end{equation}
 with $f_{x_j}$ being the partial derivative of the polynomial $f$ with respect to $x_j$. Moreover, $KR(f)$ is the module of {\it Koszul relations} spanned by obvious relations of the type $f_{x_j}f_{x_i}+(-f_{x_i})f_{x_j}=0$ and the quotient
\begin{equation} 
\label{eqER}
ER(f)=AR(f)/KR(f)
\end{equation}
is the module of {\it essential relations} (which is of course nothing else but the $n$-th cohomology group of the Koszul complex of $f_{x_0},...,f_{x_n}$), see \cite{DBull}, \cite{DSt}. Note also that with this notation, the ideal $I$ is just the saturation of the Jacobian ideal $J_f=(f_{x_0},...,f_{x_n})$.
In particular, if $V$ is a nodal hypersurface, then $I$ is just the radical of $J_f$ and hence $h \in I$ if and only if $h$ vanishes at all the nodes in $Y$, see \cite{DBull}.

It is easy to compute the dimension of the spaces $AR(f)_m$, $KR(f)_m$ and $ER(f)_m$ in terms of the dimensions of the spaces $M(f)_s$, where $M(f)=S/J_f$ is the Milnor algebra of the degree $N$ hypersurface $V:f=0$. Indeed one has the following easy result, whose proof is left to the reader. For the last claim one may see also formula (2.17) in 
\cite{DSt}.
\begin{prop}
\label{prop1} Let $g \in S_N$ be a generic polynomial (such that $g=0$ is a smooth hypersurface). Then one has the following.

\begin{enumerate}

 \item  $ \dim AR(f)_m=(n+1){n+m \choose n} -{n+m+N-1 \choose n}+ \dim M(f)_{m+N-1}.$

\item  $ \dim KR(f)_m=(n+1){n+m \choose n} -{n+m+N-1 \choose n}+ \dim M(g)_{m+N-1}.$ In particular, this number depends only on $n$ and $N$.

\item  $ \dim ER(f)_m= \dim M(f)_{m+N-1} - \dim M(g)_{m+N-1}.$

\end{enumerate}

\end{prop}

The  main result of this note is the following characterization of the fact that the singularities of $V$ are all weighted homogeneous.

\begin{thm}
\label{thm2}
Assume  that the coordinates $x$ have been chosen such that the hyperplane $H_0:x_0=0$ is transversal to $V$. Consider the projection on the first factor
\begin{equation} 
\label{proj}
p_0: AR(f)_m \to S_m/I_m,  \  \  \  (a_0,...,a_n) \mapsto [a_0].
\end{equation}
Then the equality
$$KR(f)_m=\ker p_0$$
holds for any integer $m$ if and only if all of the singularities of $V$ are weighted homogeneous.
\end{thm}

Example \ref{exCayley1} below shows the necessity of the transversality assumption.

It is clear that one has $KR(f) \subset AR(f) \cap I\Omega^n$ for any polynomial $f$. The first consequence of Theorem \ref{thm2} is a new proof of the following known result.

\begin{cor}
\label{thm1}
If all singularities of $V$ are weighted homogeneous, then
$$KR(f) = AR(f) \cap I\Omega^n.$$

\end{cor}

If we write an $n$-form $\omega \in \Omega^n$ as in the equation \eqref{eq2}, then clearly
$\omega \in I\Omega^n$ if and only if one has $\omega_K(x) \in I$ for any multi-index $K$, which says, by definition, that the polynomial $\omega_K(x)$ vanishes on the subscheme $Y$, see
 \cite{CS}, \cite{DBull}. Hence the above result says that the Koszul syzygies are in this case exactly the syzygies vanishing on the singular locus subscheme $Y$.

\begin{rk}
\label{rkA}
(i)  In the case $n=2$,  D. Cox and H. Schenck  have proven in \cite{CS} a more general version of Corollary \ref{thm1}, i.e. the partial derivatives of $f$ are replaced by a general almost complete intersection $(f_0,f_1,f_2)$. They also show that this equality holds if and only if the ideal $(f_0,f_1,f_2)$ is a local complete intersection. For a very useful discussion of the case $n=2$ we refer to \cite{ST}, subsection (1.2).

When $n>2$, a similar generalization of  Corollary \ref{thm1}  was shortly after  the paper \cite{CS} obtained by L. Bus\'e and J.-P. Jouanolou in \cite{BJ}, Proposition 4.14, in fact in a much more general setting. 

On the other hand, our Theorem \ref{thm2} can be considered to be a stronger result in this special Jacobian setting, since the vanishing on the subscheme $Y$ of just one (generic) component $a_0$ of a syzygy implies that the syzygy is in $KR(f)$,
and hence the vanising of all the other components $a_1,...,a_n$ on the subscheme $Y$.

(ii) The converse implication in Corollary \ref{thm1} seems to be an open question when $n>2$.

\end{rk}

\begin{cor}
\label{corA}
If all the singularities of $V$ are isolated and weighted homogeneous and if the hyperplane 
$H_0:x_0=0$ is transversal to $V$, then for any $m$ there is an induced monomorphism
$$\tilde p_0:ER(f)_m \to S_m/I_m.$$

\end{cor}
The equality $\dim ER(f)_m= \mu(V)$ holds as soon as $m \geq n(N-2)$, see Theorem 1 in 
\cite{DBull} and hence in this range $\tilde p_0$ is an isomorphism as $\dim S_m/I_m \leq \mu(V)$ for any $m$. Here $\mu(V)$, the global Milnor number of $V$,  is the sum of the Milnor numbers of all the singularities of $V$.
\begin{cor}
\label{corB}
If all the singularities of $V$ are isolated and weighted homogeneous, then the following inequalities hold for any integer $m$. 

\begin{enumerate}

\item $\dim ER(f)_m + \dim ER(f)_{nN-2n-1-m} \leq \mu(V)$. In particular $$\dim ER(f)_m 
\leq \mu(V)/2$$ for $m \leq (nN-2n-1)/2$. In other words the number of independent syzygies of small degree cannot be too large.

\item $\defect _mY + \defect_{nN-2n-1-m}Y \leq \mu(V)$, where $\defect_mY=\mu(V)-\dim S_m/I_m$.

\item $\dim M(f)_m + \dim M(f)_{T-m} - 2 \dim M(g)_m \leq \mu(V)$, where $T=(n+1)(N-2)$.

\end{enumerate}

\end{cor}

The first  claim that  the number of independent syzygies of small degree cannot be too large  complements the results in \cite{DS2},   \cite{Kl} showing that there are no Jacobian syzygies of 'very small' degree, see Example \ref{curves} for a precise statement in the case of curves with nodes, cusps and ordinary triple points.

The last inequality above gives a positive answer to the open Question 1 in \cite{DS1} with the additional hypothesis on the singularities of $V$ to be weighted homogeneous (equivalently, $\mu (V)=\tau (V)$, the global Tjurina number of $V$). To see this, one has to use Corollary 3 in \cite{DS1} and pay attention to the difference in notation and grading conventions between  \cite{DS1} and the present paper.

 Our interest in this type of results comes from the fact that hypersurfaces with isolated weighted homogeneous singularities have been recently shown to possess deep properties with respect to the syzygies of their Jacobian ideals, see \cite{DS1}, \cite{DS2},  \cite{DS14},  \cite{Kl}. In order to extend such properties to the non weighted homogeneous singularities, one has perhaps to understand first better the relations between syzygies and weighted homogenenous singularities.

Finally, as explained in the last section, there is a computational by-product of our result. When using a computer software as {\it Singular} \cite{Sing}, {\it CoCoA} \cite{Co} or {\it Macaulay2} \cite{Mac} to compute a basis for the vector space $AR(f)_m$ for a given $m$, usually the answer is not in the form of a basis for the subspace $KR(f)_m$ extended by some elements, whose classes yield a basis of $ER(f)_m$. Using our result, we can easily modify the basis given by the software to get a new basis satisfying the above requirement, see Examples \ref{exCayley1}, \ref{exCayley2}. Having a basis of $KR(f)_m$ is useful, since these vector spaces enter into some spectral sequences computing the topology of the hypersurface $V:f=0$ and of the associated Milnor fiber $F:f=1$, see \cite{DS1}, \cite{DS2},  \cite{DSt}.

\section{Proof of Theorem \ref{thm2} and of Corollaries  \ref{thm1}, \ref{corB}} \label{sec3}

Let the coordinates on $\PP^n$ be chosen such that the hyperplane $H_0: x_0=0$ is transversal to the hypersurface $V$ or, equivalently, that the intersection $H_0 \cap V$ is smooth. To simplify the notation, we write $f_j$ for the partial derivative $f_{x_j}$.
Then $\Gamma=V(f_1,...,f_n)$ is a 0-dimensional complete intersection contained in the affine space $U_0=\C^n=\PP^n \setminus H_0$. Moreover the ideal $J_0=(f_1,...,f_n)$ is saturated,
as it follows from Theorem 8 (Lasker's Unmixedness Theorem) in \cite{EGH}.
Then it is clear that a relation $(a_0,...,a_n) \in AR(f)$ is in $KR(f)$ if and only if $a_0 \in J_0$.
Since the ideal $J_0$ is saturated, to test the condition $a_0 \in J_0$ we can work locally, checking whether $a_{0,q} \in \tilde J_{0,q}$ at germ level, for all (closed) points $q \in \Gamma$, see the beginning of section 2 in \cite{DBull}  if necessary. Here $\tilde J_0$ denotes the sheaf ideal associated to $J_0$.

 If we use the coordinates $y_1=x_1,....,y_n=x_n$ on $U_0$, then the intersection $V_0=V \cap U_0$ is given by the equation $g(y)=0$, where
$g(y)=f(1,y_1,...,y_n)$. The Euler relation for $f$ yields the relation
\begin{equation} 
\label{euler}
f_0(1,y)+y_1g_1(y)+...+y_ng_n(y)=N \cdot g(y),
\end{equation}
where $g_j$ denotes the partial derivative of $g$ with respect to $y_j$.
If $0 \in \C^n$ is an isolated singularity of $V_0$, let $\OO_n$ be the  local ring of analytic germs at the origin, and $J_g$ the ideal in $\OO_n$ spanned by the partial derivatives of $g$. With this notation, one clearly has the following, using the characterization of isolated weighted homogeneous singularities given by K. Saito in \cite{KS}.

\begin{lem}
\label{Tj}
$\I_{Y,0}=\tilde J_{0,0}$ if and only if the germ $g$ is weighted homogeneous, i.e. $g \in J_g$.
\end{lem}

The support of $\Gamma$ consists of a finite set of points in $\PP^n$, say $q_1,...,q_r$. A part of these points, say $q_j$ for $j=1,...,s$ are  the singularities of $V$, i.e. the points in the support of $Y$.

Since $q_j$ is not a singularity for $V$ for $j>s$, it follows that $f_0(q_j)\ne 0$ in this range.
Hence, for $j>q$, the relation $R$ in \eqref{eqR} implies that the germ of function $a_{0,q_j}$ induced by $a_0$ at $q_j$ (dividing by some homogeneous polynomial $b_j$ of degree $m=\deg a_0$ such that $b_j(q_j)\ne 0$) belongs to the ideal   $\tilde J_{0,q_j}$.

Assume now that all the singularities of $V$ are weighted homogeneous, i.e. $g \in J_g$ in local coordinates at any singular point. Then
at a singular point $q_j$ with $j \leq s$, we get that $a_{0,q_j}$ belongs to the ideal  $\tilde J_{0,q_j}$ if and only if it belongs to the ideal  $\tilde I_{q_j}=\I_{Y,q_j}$, as these two stalk ideals coincide in view of Lemma  \ref{Tj}. We have shown in this way that $KR(f)_m=\ker p_0$.

Conversely, assume that some of the singularities of $V$ are not weighted homogeneous. We define a subscheme $\Gamma'$ of the scheme $\Gamma$ by describing at each point $q$ in the support of $\Gamma$ an ideal in $\OO_{\Gamma,q}$. If $q=q_j$ for $j>s$, then we set this ideal to be zero, i.e. at such a point the two schemes $\Gamma'$ and $\Gamma$ coincide.

Let now $q=g_j$ for $j \leq s$. Then $q$ is a singular point of $V$ and we may suppose that $q=0$ in order to use the above notations. Then one has $\OO_{\Gamma,0}=\OO_n/J_g$. Inside this artinian ring, there are two obvious ideals: the principal ideal spanned by the class of $g$ and the annihilator of this class. Call $\I_{\Gamma',0}$ the intersection of these two ideals and note that this ideal is nonzero exactly when the singularity $(V,0)$ is not weighted homogeneous, i.e. when $g \notin J_g$.

In this way we get a proper subscheme $\Gamma'$ of the scheme $\Gamma$ and passing to the associated saturated ideals in $S$, we get an ideal $J_0'$ containing $J_0$ and different from it.

Then, by checking locally at each point $q_j$, we see that if we choose a homogeneous polynomial $a_0 \in (J_0' \setminus J_0)$, the product $a_0f_0$ belongs to $J_0$ (this is where the annihilators come into play). Hence $a_0$ extends to a relation $(a_0,...,a_n)$ in $AR(f)$. Moreover, $p_0(a_0,...,a_n)=0$ (this is where the principal ideals spanned by the various $g$'s come into play). But by its construction, $(a_0,...,a_n) \notin KR(f)$, since $a_0 \notin J_0$. This shows that in this case the inclusion $KR(f)_m \subset \ker p_0$ is strict, and this completes the proof of Theorem \ref{thm2}.

\bigskip

We explain now how to prove Corollary \ref{thm1} using Theorem \ref{thm2}. To do this, chose a generic coordinate system $x$ such that all the coordinate hyperplanes $H_j:x_j=0$ are transversal to $V$. Then one can apply Theorem \ref{thm2} to any of the projections
\begin{equation} 
\label{projj}
p_j: AR(f)_m \to S_m/I_m,  \  \  \  (a_0,...,a_n) \mapsto [a_j]
\end{equation}
for $j=0,...,n$. It follows that $KR(f)_m=\ker p_0=...=\ker p_n$, which says exactly that 
$KR(f)_m=AR(f)_m \cap I\Omega^n$, which ends the proof of Corollary \ref{thm1}.

\bigskip

Now we give the proof of Corollary \ref{corB}. By Corollary \ref{corA} we have
$$\dim ER(f)_m \leq \dim S_m/I_m=\mu(V)-\defect_mY.$$
On the other hand, Theorem 1 in \cite{DBull} yields
$$\defect_mY=\dim ER(f)_{nN-2n-1-m}.$$
This completes the proof of the first part of the claim $(1)$ and shows that this and $(2)$ are equivalent. To prove the second part of the claim $(1)$, it is enough to recall that 
$$\dim ER(f)_m \leq \dim ER(f)_{nN-2n-1-m}$$
for $m \leq nN-2n-1-m$ in view of Corollary 11 in \cite{CD}.

To prove the third claim, we use $(1)$ and the formula given in Proposition \ref{prop1}, $(3)$.
This implies that 
$$A:=\dim M(f)_{m+N-1}-\dim M(g)_{m+N-1} +\dim M(f)_{T-m}-\dim M(g)_{T-m} \leq \mu(V),$$
where $T=(n+1)(N-2)$.
Let us denote by $B$ the RHS of the inequality $(3)$ that we want to prove.
Then, using  the equality $\dim M(g)_{T-m}=\dim M(g)_{m}$, one has
$$B=A+ (\dim M(f)_{m}-\dim M(g)_{m}) -(\dim M(f)_{m+N-1}-\dim M(g)_{m+N-1})=$$
$$=A+\dim ER(f)_{m-N+1}-\dim ER(f)_m \leq A \leq \mu(V)$$
since $\dim ER(f)_{m-N+1}\leq \dim ER(f)_{m}$ by Corollary 11 in \cite{CD}.

\section{Examples and discussion of Singular, Cocoa  and Macaulay2 outputs} \label{sec4}

\begin{ex}
\label{curves} 
Let $V:f=0$ be a plane curve  of degree $N$ with nodes $A_1$, cusps $A_2$ and ordinary triple points $D_4$. Then it follows from Example 2.2 in \cite{DS14} that one has
$$AR(f)_m=ER(f)_m=0$$
 for $m<2N/3-2$. Corollary \ref{corB}, $(1)$ implies that we also have
$$\dim AR(f)_m=\dim ER(f)_m \leq \frac{n(A_1)+2n(A_2)+4n(D_4)}{2}$$
for $m \leq N-3$, where $n(\Sigma)$ denotes the number of singularities of the curve $V$ of type $\Sigma$. For example the line arrangement given by $f=(x^2-y^2)(x^2-z^2)(y^2-z^2)=0$ has 3 nodes and 4 triple points and the above inequality becomes 
$$\dim AR(f)_3=\dim ER(f)_3 \leq 9.$$
A computation using for instance {\it Singular} gives $\dim AR(f)_3=\dim ER(f)_3 =4$.

\end{ex}

\begin{ex}
\label{exCayley1} Consider the Cayley cubic surface $V$ in $\PP^3$ with coordinates $(x,y,z,w)$ given by the equation $f=xyz+xyw+xzw+yzw=0$. This surface has four $A_1$ singularities located at $q_1=(1:0:0:0)$, $q_2=(0:1:0:0)$, $q_3=(0:0:1:0)$ and $q_4=(0:0:0:1)$

A direct computation using {\it Singular } shows that the Hilbert series of the graded Milnor algebras $M(g)$ and $M(f)$ in this case are
$$H(M(g),t)=1+4t+6t^{2}+4t^{3}+t^{4},$$
and
$$H(M(f),t)=1+4t+6t^{2}+4t^{3}+4t^4+ ...$$
Using the  formulas given in Proposition \ref{prop1}, with $n=N=3$ and $m=2$, we get $\dim AR(f)_2=9$, $\dim KR(f)_2=6$ and $\dim ER(f)_2=3$. We also get $\dim ER(f)_k=4=\mu(V)$ for any $k \geq 3$.

First we illustrate Corollary \ref{thm1}.
 The {\it Singular}  package gives the following list of syzygies of degree two
among the partial derivatives $f_x, f_y, f_z, f_w$, i.e. a basis for $AR(f)_2$.
$$s_1=
(-2xw-zw,2yz+zw,-2z^2-zw, zw+2w^2),$$
$$s_2=
(-2yz+2xw+2yw+zw,-2yw-zw,-2yz-2z^2-2yw-3zw,2yz+2yw+3zw+2w^2),$$
$$s_3=
(xz-xw,yz-yw,-2z^2-zw,zw+2w^2),$$
$$s_4=
(-2xw-zw,-2xz+2xw+2yw+zw,-2xz-2z^2-2xw-3zw,2xz+2xw+3zw+2w^2),$$
$$s_5=
(2xw+yw,2y^2+yw,-2yz-yw,-yw-2w^2),$$
$$s_6=
(2xy+4xw+3yw,2y^2+yw,-4yz-3yw-2zw,-yw-2w^2),$$
$$s_7=
(0,xy+xw+yw,-xz-xw-zw,0),$$
$$s_8=
(0,0,-xy-xz-yz,xy+xw+yw),$$
$$s_9=
(2x^2+xw,2xy+3xw+4yw,-4xz-3xw-2zw,-xw-2w^2).$$
With obvious meaning, we have the equalities
$$s_j(q_1)=(0,0,0,0), \\ s_j(q_2)=(0,0,0,0), \\ s_j(q_3)=(0,0,-2,0), \\ s_j(q_4)=(0,0,0,2)$$
for $j=1,2,3,4$,
$$s_k(q_1)=(0,0,0,0), \\ s_k(q_2)=(0,2,0,0), \\ s_k(q_3)=(0,0,0,0), \\ s_k(q_4)=(0,0,0,-2)$$
for $k=5,6$,
$$s_m(q_1)=(0,0,0,0), \\ s_m(q_2)=(0,0,0,0), \\ s_m(q_3)=(0,0,0,0), \\ s_m(q_4)=(0,0,0,0)$$
for $m=7,8$ and finally
$$s_9(q_1)=(2,0,0,0), \\ s_9(q_2)=(0,0,0,0), \\ s_9(q_3)=(0,0,0,0), \\ s_9(q_4)=(0,0,0,-2).$$
It follows that $s_1-s_2$, $s_1-s_3$, $s_1-s_4$, $s_5-s_6$, $s_7$ and $s_8$ form a basis for the vector space $KR(f)_2$, while the classes of $s_1$, $s_5$ and $s_9$ form a basis for $ER(f)_2$. Note that the coordinate hyperplane $x=0$ is not transversal to $V$ in this case and the morphism $\tilde p_0$ from Theorem \ref{thm2} is not injective, since for instance $\tilde p_0(s_1)=0.$
\end{ex}

\begin{ex}
\label{exCayley2} 

Now we illustrate Theorem \ref{thm2}.
Consider  the same Cayley cubic surface $V$ in $\PP^3$ (obtained from the previous one by a linear change of coordinates)  with coordinates denoted again $(x,y,z,w)$ for simplicity, having four $A_1$ singularities located this time at $q_1=(-1:1:1:1)$, $q_2=(1:-1:1:1)$, $q_3=(1:1:-1:1)$ and $q_4=(1:1:1:-1)$ and given by
$f=4(x^3+y^3+z^3+w^3)-(x+y+z+w)^3$. Now the coordinate hyperplane $x=0$ is  transversal to $V$.

The {\it Singular} package gives the  list of syzygies $s_1,...,s_9$ of degree two
among the partial derivatives $f_x, f_y, f_z, f_w$, i.e. a basis for $AR(f)_2$, and we list below only the polynomial $t_j$, the corresponding coefficient $a_0$ in the syzygy $s_j$, for $j=1,...,9$.
$$t_1=
-yz+z^2+xw+2yw+2zw+w^2,$$
$$t_2=
2xz+3yz+3z^2+xw+6zw+3w^2,$$
$$t_3=0$$
$$t_4=
2xy+6y^2+2xz+6yz+3z^2+2xw+6yw+6zw+3w^2,$$
$$t_5=
-y^2+yz-xw-2yw-2zw-w^2,$$
$$t_6=
2xy+3y^2+3yz+xw+6yw+3w^2,$$
$$t_7=
3x^2+8xy-3y^2+8xz+12yz+6z^2+8xw+12yw+12zw+6w^2,$$
$$t_8=
y^2-z^2,$$
$$t_9=
4x^2+8xy+3y^2+8xz+5yz+4z^2+11xw+14yw+14zw-5w^2.$$
We have the equalities
$$(t_1(q_1),t_1(q_2),t_1(q_3),t_1(q_4))     =(4,4,4,-4),$$
$$(t_2(q_1),t_2(q_2),t_2(q_3),t_2(q_4))     =(12,12,-4,4),$$
$$(t_3(q_1),t_3(q_2),t_3(q_3),t_3(q_4)) =(t_8(q_1),t_8(q_2),t_8(q_3),t_8(q_4))     =(0,0,0,0),$$
$$(t_4(q_1),t_4(q_2),t_4(q_3),t_4(q_4))     =(24,8,8,8),$$
$$(t_5(q_1),t_5(q_2),t_5(q_3),t_5(q_4))     =(-4,-4,-4,4),$$
$$(t_6(q_1),t_6(q_2),t_6(q_3),t_6(q_4))     =(12,-4,12,4),$$
$$(t_7(q_1),t_7(q_2),t_7(q_3),t_7(q_4))     =(24,8,8,8),$$
$$(t_9(q_1),t_9(q_2),t_9(q_3),t_9(q_4))     =(12,12,12,-12).$$
It follows that $s_3$, $s_5+s_1$, $s_6-s_4+s_2$, $s_7-s_4$, $s_8$ and $s_9-3s_1$
 form a basis for the vector space $KR(f)_2$, while the classes of $s_1$, $s_2$ and $s_4$ form a basis for $ER(f)_2$.

The classes of $xs_1$, $xs_2$ and $xs_4$ are linearly independent in $ER(f)_3$ by Corollary 11 in
\cite{CD}.  Using Theorem \ref{thm2} we see that by adding the class of $ws_1$ we get a basis of 
$ER(f)_3$. Applying again Corollary 11, we see that the classes of $x^{k-2}s_1$, $x^{k-2}s_2$, $x^{k-2}s_4$ and $x^{k-3}ws_1$ form a basis for $ER(f)_k$ for any $k \geq 3$.
\end{ex}

\begin{rk}
\label{rkC}
In the above examples we have discussed the outputs given by the software {\it Singular}. Completely similar features are shared by the outputs produced by {\it CoCoA} and {\it Macaulay2}. We can provide the corresponding listings on request.

\end{rk}

\end{document}